\documentclass[12pt]{article}

\usepackage{amsmath}
\usepackage{theorem}
\usepackage{amssymb}
\usepackage{latexsym}
\usepackage{a4wide}
\usepackage{longtable}
\usepackage{epic}

\newtheorem{thm}{Theorem}

{\theorembodyfont{\rmfamily} }
{\theorembodyfont{\rmfamily} \newtheorem{exa}[thm]{Example}}

\newcommand{\C}{\mathbb{C}}

\newcommand{\Z}{\mathbb{Z}}

\newcommand{\ad}{\mathrm{\mathop{ad}}}

\newcommand{\ssl}{\mathfrak{\mathop{sl}}}
\newcommand{\g}{\mathfrak{g}}

\newcommand{\lab}[1]{\footnotesize{#1}}

\begin{document}

\title{Computations with reachable elements in simple Lie algebras}
\author{Willem A. de Graaf\\
Dipartimento di Matematica\\
Universit\`{a} di Trento\\
Italy}
\date{}
\maketitle

\begin{abstract}
We report on some computations with reachable elements in simple Lie
algebras of exceptional type within the {\sf SLA} package of {\sf GAP}4. 
These computations confirm the classification of such elements by
Elashvili and Gr\'elaud. Secondly they answer a question from Panyushev.
Thirdly they show in what way a recent result of Yakimova for the
Lie algebras of classical type extends to the exceptional types.  
\end{abstract}

\section{Introduction}

Let $\g$ be a simple Lie algebra over $\C$ (or over an algebraically
closed field of characteristic 0). For $e\in \g$ we denote its centraliser
in $\g$ by $\g_e$. In \cite{panyushev4} an $e$ in $\g$ is defined to be
{\em reachable} if $e \in [\g_e,\g_e]$. Such an element has to be nilpotent.
In \cite{elashgre}, Elashvili and Gr\'elaud gave a classification of reachable 
elements in $\g$ (in this paper such elements are called {\em compact}, 
in analogy with \cite{blabry}).  

By the Jacobson-Morozov theorem a nilpotent $e\in \g$ lies in an
$\ssl_2$-triple $(h,e,f)$ (where $[e,f]=h$, $[h,e]=2e$, $[h,f]=-2f$).
By the adjoint representation the subalgebra spanned by such a triple acts
on $\g$. Since the eigenvalues of $\ad h$ are integers, we get a grading
$$ \g = \bigoplus_{k\in\Z} \g(k)$$
where $\g(k) = \{ x \in \g\mid [h,x]=kx\}$. Now set $\g(k)_e = \g(k)\cap 
\g_e$, and let $\g(\geq 1)_e$ denote the subalgebra spanned by all
$\g(k)_e$, $k\geq 1$. 

Panyushev (\cite{panyushev4}) showed that, for $\g$ of type $A_n$,
$e$ is reachable if and only if $\g(\geq 1)_e$ is generated as Lie algebra
by $\g(1)_e$. Here we call this the {\em Panyushev property} of $\g$.
In \cite{panyushev4} it is stated that this property also holds for the 
other classical types and the question is posed whether it holds for the
exceptional types. In \cite{yakimova2} a proof is given that the Panyushev
property holds in types $B_n$, $C_n$, $D_n$. 
Our computations confirm that the Panyushev property holds also for
the Lie algebras of exceptional type.

Yakimova (\cite{yakimova2}) studied the stronger condition $\g_e = 
[\g_e,\g_e]$. For the purposes of this paper we call elements $e$
satisfying this condition {\em strongly reachable}. She showed that
for $\g$ of classical type, $e$ is strongly reachable if and only if
the nilpotent orbit of $e$ is rigid. (This means that it is not induced,
cf. \cite{lusp}, \cite{elashvili}, \cite{elasgra}.) Furthermore, 
this is shown to fail for $\g$ of exceptional type. As a result of our
calculations we find all rigid nilpotent orbits whose representatives are
not strongly reachable. From this we conclude that $e$ is strongly 
reachable if and only if $e$ is both reachable and rigid.  We note that 
one direction of this statement can be shown in a uniform way for all 
$\g$: if $e$ is strongly reachable then it is reachable, but also rigid by
\cite{yakimova2}, Proposition 11. The converse for exceptional types follows
from our calculations in two ways. Firstly we compute the list of all
strongly reachable orbits and the list of all nilpotent orbits that are
reachable and rigid, and find that they are the same. Second, the Panyushev 
property, which we checked by computation for the exceptional types, 
also implies the statement. For the classical types we have, of course,
the stronger theorem from \cite{yakimova2}.

The {\sf SLA} package (\cite{sla}), written in the language of the computer 
algebra system {\sf GAP}4 (\cite{gap4}), has functionality for working
with the nilpotent orbits in simple Lie algebras. In particular the package
contains the classifcation of such orbits. Using this it is 
straightforward to approach the above questions by computational means.
Indeed, for a nilpotent orbit the system easily computes a representative
$e$, and a corresponding $\ssl_2$-triple.
Then using functions present in {\sf GAP}4 we can compute the centralizer,
$\g_e$, and its derived subalgebra, and check whether $e$ lies in it.
This gives us the list of reachable nilpotent orbits. Secondly, a similar
procedure yields the list of strongly reachable orbits.
Thirdly, {\sf SLA} has
a function for computing the grading corresponding to an $\ssl_2$-triple.
With that it is straightforward to check whether $\g(\geq 1)_e$ is generated 
by $\g(1)_e$. The appendix contains the code for the functions implementing
these procedures.

{\bf Acknowledgement:} I thank Alexander Elashvili for suggesting the
topics of this paper to me.

\section{Reachable nilpotent elements in the Lie algebras of exceptional
type}

Tables \ref{tab:rigidE6}, \ref{tab:rigidE7}, \ref{tab:rigidE8}, 
\ref{tab:rigidF4}, and \ref{tab:rigidG2} contain the nilpotent orbits
that by our calculations are reachable. The content of the tables is as
follows. The first column has the label of the orbit, and the second
column the weighted Dynkin diagram. The third and fourth columns contain
a $\times$ if the orbit is, respectively, strongly reachable and rigid.
We note that the classification of rigid nilpotent orbits is known
(see \cite{elashvili}, \cite{elasgra}). 

\setlongtables

\begin{longtable}{|l|c|c|c|}
\caption{Reachable nilpotent orbits in $E_6$.}\label{tab:rigidE6}
\endfirsthead
\hline
\multicolumn{4}{|l|}{\small\slshape Reachable nilpotent orbits in $E_6$.} \\
\hline
\endhead
\hline
\endfoot
\endlastfoot

\hline

label & characteristic & Strong & Rigid\\
\hline

$A_1$ & $0~~~~0~~~~\overset{\text{\normalsize 1}}{0}~~~~0~~~~0$ 
& $\times$ & $\times$ \\

$2A_1$ & $1~~~~0~~~~\overset{\text{\normalsize 0}}{0}~~~~0~~~~1$ 
&&\\

$3A_1$ & $0~~~~0~~~~\overset{\text{\normalsize 0}}{1}~~~~0~~~~0$ 
& $\times$ & $\times$ \\

$A_2+A_1$ & $1~~~~0~~~~\overset{\text{\normalsize 1}}{0}~~~~0~~~~1$ 
&&\\

$A_2+2A1$ & $0~~~~1~~~~\overset{\text{\normalsize 0}}{0}~~~~1~~~~0$ 
&&\\

$2A_2+A_1$ & $1~~~~0~~~~\overset{\text{\normalsize 0}}{1}~~~~0~~~~1$ 
& $\times$ & $\times$ \\

\hline

\end{longtable}

\begin{longtable}{|l|c|c|c|}
\caption{Reachable nilpotent orbits in $E_7$.}\label{tab:rigidE7}
\endfirsthead
\hline
\multicolumn{4}{|l|}{\small\slshape Reachable nilpotent orbits in $E_7$.} \\
\hline
\endhead
\hline
\endfoot
\endlastfoot

\hline
label & characteristic & Stong & Rigid\\
\hline

$A_1$ & $1~~~~0~~~~\overset{\text{\normalsize 0}}{0}~~~~0~~~~0~~~~0$ 
& $\times$ & $\times$\\

$2A_1$ & $0~~~~0~~~~\overset{\text{\normalsize 0}}{0}~~~~0~~~~1~~~~0$ 
& $\times$ & $\times$\\

$(3A_1)'$ & $0~~~~1~~~~\overset{\text{\normalsize 0}}{0}~~~~0~~~~0~~~~0$ 
& $\times$ & $\times$\\

$4A_1$ & $0~~~~0~~~~\overset{\text{\normalsize 1}}{0}~~~~0~~~~0~~~~1$ 
& $\times$ & $\times$\\

$A_2+A_1$ & $1~~~~0~~~~\overset{\text{\normalsize 0}}{0}~~~~0~~~~1~~~~0$ 
&&\\

$A_2+2A_1$ & $0~~~~0~~~~\overset{\text{\normalsize 0}}{1}~~~~0~~~~0~~~~0$ 
& $\times$ & $\times$\\

$2A_2+A_1$ & $0~~~~1~~~~\overset{\text{\normalsize 0}}{0}~~~~0~~~~1~~~~0$ 
& $\times$ & $\times$\\

$A_4+A_1$ & $1~~~~0~~~~\overset{\text{\normalsize 0}}{1}~~~~0~~~~1~~~~0$ 
&&\\

\hline

\end{longtable}

\begin{longtable}{|l|c|c|c|}
\caption{Reachable nilpotent orbits in $E_8$.}\label{tab:rigidE8}
\endfirsthead
\hline
\multicolumn{4}{|c|}{\small\slshape Reachable nilpotent orbits in $E_8$.} \\
\hline
\endhead
\hline
\endfoot
\endlastfoot

\hline
label & characteristic & Stong & Rigid\\
\hline

$A_1$ &
$0~~~~0~~~~\overset{\text{\normalsize 0}}{0}~~~~0~~~~0~~~~0~~~~1$ 
& $\times$ & $\times$\\

$2A_1$ &
$1~~~~0~~~~\overset{\text{\normalsize 0}}{0}~~~~0~~~~0~~~~0~~~~0$
& $\times$ & $\times$\\

$3A_1$ &
$0~~~~0~~~~\overset{\text{\normalsize 0}}{0}~~~~0~~~~0~~~~1~~~~0$
& $\times$ & $\times$\\

$4A_1$ &
$0~~~~0~~~~\overset{\text{\normalsize 1}}{0}~~~~0~~~~0~~~~0~~~~0$
& $\times$ & $\times$\\

$A_2+A_1$ &
$1~~~~0~~~~\overset{\text{\normalsize 0}}{0}~~~~0~~~~0~~~~0~~~~1$
& $\times$ & $\times$\\

$A_2+2A_1$ &
$0~~~~0~~~~\overset{\text{\normalsize 0}}{0}~~~~0~~~~1~~~~0~~~~0$
& $\times$ & $\times$\\

$A_2+3A_1$ &
$0~~~~1~~~~\overset{\text{\normalsize 0}}{0}~~~~0~~~~0~~~~0~~~~0$
& $\times$ & $\times$\\

$2A_2+A_1$ &
$1~~~~0~~~~\overset{\text{\normalsize 0}}{0}~~~~0~~~~0~~~~1~~~~0$
& $\times$ & $\times$\\

$A_4+A_1$ &
$1~~~~0~~~~\overset{\text{\normalsize 0}}{0}~~~~0~~~~1~~~~0~~~~1$
&&\\

$2A_2+2A_1$ &
$0~~~~0~~~~\overset{\text{\normalsize 0}}{0}~~~~1~~~~0~~~~0~~~~0$
& $\times$ & $\times$\\

$(A_3+2A_1)''$ &
$0~~~~1~~~~\overset{\text{\normalsize 0}}{0}~~~~0~~~~0~~~~0~~~~1$
& $\times$ & $\times$\\

$D_4(a_1)+A_1$ &
$0~~~~0~~~~\overset{\text{\normalsize 1}}{0}~~~~0~~~~0~~~~1~~~~0$
& $\times$ & $\times$\\

$A_3+A_2+A_1$ &
$0~~~~0~~~~\overset{\text{\normalsize 0}}{1}~~~~0~~~~0~~~~0~~~~0$
& $\times$ & $\times$\\

$2A_3$ &
$1~~~~0~~~~\overset{\text{\normalsize 0}}{0}~~~~1~~~~0~~~~0~~~~0$ 
& $\times$ & $\times$\\

$A_4+2A_1$ &
$0~~~~0~~~~\overset{\text{\normalsize 0}}{1}~~~~0~~~~0~~~~0~~~~1$ 
&&\\

$A_4+A_3$ &
$0~~~~0~~~~\overset{\text{\normalsize 0}}{1}~~~~0~~~~0~~~~1~~~~0$ 
& $\times$ & $\times$\\

\hline

\end{longtable}

\begin{longtable}{|l|c|c|c|}
\caption{Reachable nilpotent orbits in $F_4$.}\label{tab:rigidF4}
\endfirsthead
\hline
\multicolumn{4}{|l|}{\small\slshape Reachable nilpotent orbits in $F_4$.} \\
\hline
\endhead
\hline
\endfoot
\endlastfoot

\hline

label & characteristic & Strong & Rigid\\
\hline
 &
\begin{picture}(20,7)
  \put(-20,0){\circle{6}}
  \put(0,0){\circle{6}}
  \put(20,0){\circle{6}}
  \put(40,0){\circle{6}}
  \put(-17,0){\line(1,0){14}}
  \put(2,2){\line(1,0){16}}
  \put(2,-2){\line(1,0){16}}
\put(5,-3){$>$}
  \put(23,0){\line(1,0){14}}
\end{picture} & &

\\
\hline

$A_1$ & 1~~~0~~~0~~~0
& $\times$ & $\times$\\

$\widetilde{A}_1$ & 0~~~0~~~0~~~1
& $\times$ & $\times$\\

$A_1+\widetilde{A}_1$ & 0~~~1~~~0~~~0
& $\times$ & $\times$\\

$A_2+\widetilde{A}_1$ & 0~~~0~~~1~~~0
& $\times$ & $\times$\\

\hline

\end{longtable}

\begin{longtable}{|l|c|c|c|}
\caption{Reachable nilpotent orbits in $G_2$.}\label{tab:rigidG2}
\endfirsthead
\hline
\multicolumn{4}{|l|}{\small\slshape Reachable nilpotent orbits in $G_2$.} \\
\hline
\endhead
\hline
\endfoot
\endlastfoot

\hline

label & characteristic & Strong & Rigid\\

&

  \begin{picture}(10,7)
  \put(-10,0){\circle{6}}
  \put(22,0){\circle{6}}
\put(3,-3){$>$}
  \put(-8,-2){\line(1,0){28}}
  \put(-7,0){\line(1,0){26}}
  \put(-8,2){\line(1,0){28}}
\end{picture} &&                   \\
\hline

\hline

$\widetilde{A}_1$ & 1~~0 
& $\times$ & $\times$
\\

\hline

\end{longtable}

We make the following comments.

\begin{itemize}
\item Here the reachable elements are exactly the same as in the
paper of Elashvili and Gr\'elaud. Therefore our calculations confirm their
result.
\item The rigid nilpotent orbits that are not strongly reachable are
\begin{itemize}
\item in type $E_7$: $(A_3+A_1)'$ $(41,40)$,
\item in type $E_8$: $A_3+A_1$ $(84,83)$, $D_5(a_1)+A_2$ $(46,45)$, 
$A_5+A_1$ $(46,45)$,
\item in type $F_4$: $\widetilde{A}_2+A_1$ $(16,15)$,
\item in type $G_2$: $A_1$ $(6,5)$.
\end{itemize}
Here the pair of integers in brackets is $(\dim \g_e,\dim [\g_e,\g_e])$.

\item In type $E_6$ all rigid orbits are strongly reachable. Hence 
in this type the situation is the same as for the classical types: 
$e$ is strongly reachable if
and only if the orbit of $e$ is rigid.
\item The last two columns of all tables are equal. This shows that
for the exceptional types the following theorem holds: $e$ is strongly
reachable if and only if $e$ is both reachable and rigid. 
\item This last statement also follows from the Panyushev property. 
Indeed, $e$ rigid implies that $\g(0)_e$ is semisimple, so
$[\g(0)_e,\g(0)_e]=\g(0)_e$. Furthermore, $[\g(0)_e,\g(1)_e]= 
\g(1)_e$ by \cite{yakimova2}, Lemma 8 (where this is shown to hold for 
all nilpotent $e$).
By the Panyushev property this implies that $[\g_e,\g_e]= \g_e$.
\item We see that for all nilpotent orbits that are rigid but not
stronly reachable the codimension of $[\g_e,\g_e]$ in $\g_e$ is 1.
Since a rigid orbit is reachable if and only if it is strongly reachable,
we get that in all those cases $e$ spans the quotient $\g_e/[\g_e,\g_e]$.
\end{itemize}

\begin{exa}
Let us consider the nilpotent orbit in the Lie algebra of type $E_7$ with
label $A_3+A_2$. This orbit is not reachable. It has a representative
with diagram

\begin{picture}(5,35) 
\put(3,10){\circle{6}} 
\put(5,15){\lab{29}} 
\put(23,10){\circle{6}}
\put(25,15){\lab{32}}
\put(6,10){\line(1,0){14}}
\put(43,10){\circle{6}}
\put(26,10){\line(1,0){14}}
\put(45,15){\lab{31}}
\put(63,10){\circle{6}}
\put(65,15){\lab{27}}
\put(83,10){\circle{6}}
\put(85,15){\lab{30}}
\put(66,10){\line(1,0){14}}
\end{picture}

This means that the representative is $e=x_{29}+x_{32}+x_{31}+x_{27}+x_{30}$,
where $x_i$ denotes the root vector corresponding to the $i$-th positive 
root (enumeration as in {\sf GAP}4, cf. \cite{gra14}). Furthermore,
the Dynkin diagram of these roots is as shown above. This representative
is stored in the package {\sf SLA}.

Now, if the orbit were reachable then $e \in [\g_e,\g_e]\cap \g(2)$.
Using the {\sf SLA} package we can easiliy compute the latter space:

\begin{verbatim}
gap> L:= SimpleLieAlgebra("E",7,Rationals);;
gap> o:= NilpotentOrbits(L);;
gap> sl2:=SL2Triple( o[19] );
[ (2)*v.90+(3)*v.92+(2)*v.93+(3)*v.94+(4)*v.95, (6)*v.127+(9)*v.128+(12)*v.129
+(18)*v.130+(14)*v.131+(10)*v.132+(5)*v.133, v.27+v.29+v.30+v.31+v.32 ]
gap> g:= SL2Grading( L, sl2[2] );;
gap> g2:= Subspace( L, g[1][2] );;
gap> der:= LieDerivedSubalgebra(LieCentralizer(L,Subalgebra(L,[sl2[3]])));
<Lie algebra of dimension 33 over Rationals>
gap> BasisVectors( Basis( Intersection( g2, der ) ) );
[ v.18, v.23+(-1)*v.24+v.28, v.24+(-1)*v.25+(-1)*v.28, 
v.27+(-1)*v.29+v.30+(-1)*v.31+(-1)*v.32, v.33+(-1)*v.36+v.37, 
v.34+(-1)*v.36+v.37, v.39 ]
\end{verbatim}

First we make some comments on the above computation. The $\ssl_2$-triple
comes ordered as $(f,h,e)$. So the second element is the neutral element,
and the third element is the nil-positive element, i.e., the representative,
which is as indicated above. So the second element defines the grading,
which we compute with {\tt SL2Grading}. In the subsequent line the
subspace $\g(2)$ is defined, followed by $[\g_e,\g_e]$. Finally a basis
of the intersection is computed.

We see that one of the basis vectors of the intersection is 

\begin{verbatim}
v.27+(-1)*v.29+v.30+(-1)*v.31+(-1)*v.32.
\end{verbatim}

So we see that $e=(x_{29}+x_{32}+x_{31})+(x_{27}+x_{30})$ does not lie
in $[\g_e,\g_e]$ but $(x_{29}+x_{32}+x_{31})-(x_{27}+x_{30})$ does!
\end{exa}

\section*{Appendix: the code}

\begin{verbatim}
ReachableOrbits:= function( L )

    # this returns the nilpotent orbits of L that are reachable.

    local o, reachables, i, sl2, e, K;

    o:= NilpotentOrbits(L);
    reachables:= [ ];
    for i in [1..Length(o)] do
       sl2:= SL2Triple( o[i] );
       e:= sl2[3];
       K:= LieCentralizer( L, Subalgebra(L,[e]) );
       if e in LieDerivedSubalgebra(K) then
          Add( reachables, o[i] );
       fi;
    od;

    return reachables;

end;

PanyushevProperty:= function( L )

    # this function returns true if the Panyushev property 
    # holds for L, otherwise false is returned. 

    local reachables, sl2, K, prop, r, c, M, g;
    
    reachables:= ReachableOrbits(L);
    prop:= true;
    for r in reachables do
        sl2:= SL2Triple( r );
        g:= SL2Grading( L, sl2[2] );
        K:= LieCentralizer( L, Subalgebra(L,[sl2[3]]) );
        c:= List( g[1], u -> BasisVectors( Basis( Intersection( 
                            K, Subspace(L,u) ) ) ) );
        M:= Subalgebra( L, Flat( c ) ); 
        if Dimension( Subalgebra( L, c[1] ) ) <> Dimension(M) then
           Print("Property not verified for  ",r[1][3],"\n");
           prop:= false;
        fi;
    od;

    return prop;

end;


StronglyReachableOrbits:= function( L )

    # returns the stronly reachable nilpotent orbits of L.
   
    local o, good, i, sl2, e, K;

    o:= NilpotentOrbits(L);
    good:= [ ];
    for i in [1..Length(o)] do
       sl2:= SL2Triple( o[i] );
       e:= sl2[3];
       K:= LieCentralizer( L, Subalgebra(L,[e]) );
       if LieDerivedSubalgebra(K)=K then
          Add( good, o[i] );
       fi;
    od;

    return good;

end;
\end{verbatim}

As a small illustration we give a sample session for the Lie algebra
of type $E_6$.

\begin{verbatim}
gap> RequirePackage("sla");
gap> L:= SimpleLieAlgebra("E",6,Rationals);;
gap> r:= ReachableOrbits( L );;
gap> List( r, WeightedDynkinDiagram );
[ [ 0, 1, 0, 0, 0, 0 ], [ 1, 0, 0, 0, 0, 1 ], [ 0, 0, 0, 1, 0, 0 ], 
  [ 1, 1, 0, 0, 0, 1 ], [ 0, 0, 1, 0, 1, 0 ], [ 1, 0, 0, 1, 0, 1 ] ]
gap> r:= StronglyReachableOrbits( L );;
gap> List( r, WeightedDynkinDiagram );
[ [ 0, 1, 0, 0, 0, 0 ], [ 0, 0, 0, 1, 0, 0 ], [ 1, 0, 0, 1, 0, 1 ] ]
gap> PanyushevProperty(L);
true
\end{verbatim}

\def\cprime{$'$} \def\cprime{$'$} \def\Dbar{\leavevmode\lower.6ex\hbox to
  0pt{\hskip-.23ex \accent"16\hss}D} \def\cprime{$'$} \def\cprime{$'$}
  \def\cprime{$'$}

\end{document}